\documentclass[12pt]{article}
\usepackage{amssymb}
\usepackage{tikz-cd}
\usepackage{amscd}
  \usepackage{graphicx}
 \usepackage{amsmath,amsthm,amssymb,amsfonts}
  \usepackage{indentfirst}
\newtheorem{theorem}{Theorem}

\newtheorem{proposition}{Proposition}

\newtheorem{remark}{Remark}
\newtheorem{definition}{Definition}
\newtheorem{lemma}{Lemma}

\newcommand{\CC}{{\mathbb C}}

\newcommand{\Z}{{\mathbb Z}}
\newcommand{\R}{{\mathbb R}}

\newcommand{\CP}{{\mathbb {CP}}}

\title{Picard-Lefschetz Monodromy Groups of Quadratic
Hypersurfaces}
\author{Daodao Yang}
\begin{document}
\date{}
\maketitle

\centerline{\bf Abstract}
\medskip

We study the topology of the space of affine hyperplanes $L \subset \CC^n$
which are in general position with respect to a given generic quadratic
hypersurface $A$, and calculate the monodromy action of the fundamental group
of this space on the relative homology groups $H_*(\CC^n, A \cup L)$ associated
with such hyperplanes. \bigskip

\section{The statement of the problem and the relative homology group}

$ A $ is an non-degenerate quadratic hypersurface in $ \CC^{n} $.
\medskip

For instance, $ A $  could be the set $ \lbrace (z_{1}, z_{2}, ... z_{n}) \in
\CC^{n}~\Big |~   z_{1}^{2}+ z_{2}^{2}+ ... +z_{n}^{2} = 1\rbrace $.\medskip

$ L $ is a complex hyperplane in $\CC^n$.

By $\CP^{n-1}_{\infty}$ we denote the ``infinitely distant'' part $\CP^{n}
\setminus \CC^{n}$ of the projective closure of $\CC^n$.

$ \overline{A} $ is the closure of $ A $ in $\CP^{n}$. Non-degeneracy of $A$
implies that $\overline{A}$ is smooth in $\CP^n$ and intersects $\CP^{n-1}_\infty$ transversally, and so
$\overline{A} \cap \CP^{n-1}_\infty$ is a non-degenerate quadric hypersurface
in $\CP^{n-1}$.

Let $\overset{\vee}{\CP^{n}} $ be the space of all hyperplanes in $ \CP^{n}$.

\begin{definition} \rm
$L$ is asymptotic for $ A \subset \CC^{n} $ if $ \overline{L}\cap
\CP^{n-1}_{\infty} $ is tangent to $\overline{A}\cap\CP^{n-1}_{\infty}$.

$L$ is \underline{not} in general position with respect to $ A $ if either it
is tangent to $A$ at some point in $\CC^n$, or it is asymptotic for $A$.
\end{definition}

In other words, $L$ is in general position with respect to $A$ if and only if
its closure $\overline{L} \subset \CP^n$ is transversal to the (stratified) algebraic
set $A \cup \CP^{n-1}_\infty$.
\medskip

{\bf Notation.} Denote by $\overset{\vee}{A}$ and $ \overset{\ast}{A}$ the subsets in
$\overset{\vee}{\CP^{n}} $ consisting of all tangent and asymptotic hyperplanes
of $A$, respectively; in addition, the point in $\overset{\vee}{\CP^{n}} $
corresponding to the ``infinitely distant'' hyperplane also is by definition included into $
\overset{\ast}{A}$.
\bigskip

By the Thom's isotopy lemma (see \cite{GM}, \cite{APLT}) the pairs of spaces
$(\CP^n, \overline{A} \cup \overline{L} \cup \CP^{n-1}_{\infty})$ form a locally trivial fiber bundle over the
space $\overset{\vee}{\CP^{n}} \setminus (\overset{\vee}{A} \cup
\overset{\ast}{A})$ of planes $L$ which are in general position with respect to
$A$. Therefore the fundamental group of the latter space acts  on all
homology groups related with spaces $(\CP^n, \overline{A}  \cup \overline{L} \cup \CP^{n-1}_{\infty} )$ by the monodromy, in particular on the groups
$H_n(\CC^n, A \cup L)$. The explicit calculation of this action is the main
goal of this work; this is a sample result for a large family of similar problems concerning the hypersurfaces of higher degrees and/or non-generic ones.

This action is important in the problems of integral geometry, when the
integration contour is represented by a relative chain in $\CC^n$ with boundary
at $A \cup L$, and integration $n$-form is holomorphic and has singularity at
the infinity; see e.g. \cite{APLT}, Chapter III.

We always assume that $n \geq 2$, because otherwise the problem is trivial.

\subsection{The representation space}

\begin{proposition} \label{relativehomology}
If $L$ is in general position with respect to $A$, then $$ H_{n}(\CC^{n}, A\cup
L) \cong H_{n-1}(A\cup L) \cong\Z^{2}, $$ and $ H_{i}(\CC^{n}, A\cup L) \cong
\tilde H_{i-1}(A\cup L) \cong 0$ for all $i \neq n$ $($here $\tilde H$ means
homology group reduced modulo a point$)$.
\end{proposition}

Proof. First, we have the long exact sequence for the pair $ (\CC^{n}, A\cup L) $:

\begin{equation}
\label{bpr}
...\rightarrow H_{i}(A \cup L)
\rightarrow
H_{i}(\CC^{n})
\rightarrow
H_{i}(\CC^{n}, A \cup L) \rightarrow H_{i-1}(A \cup L)
\rightarrow
H_{i-1}(\CC^{n})
\rightarrow
...\end{equation}
\bigskip

The homology groups of $\CC^{n}$ coincide with these of a point. So $
H_{i}(\CC^{n}, A\cup L) \cong \tilde H_{i-1}(A\cup L)$ for any $i$.
\bigskip

Second, Milnor theorem shows that $ A $ is homotopy equivalent to $ S^{n-1} $,
and $ A \cap L $ is homotopy equivalent to $ S^{n-2} $. Thus $H_{k}(A) =
\begin{cases}0& \mbox{for } \emph{\emph{others}}; \cr \mathbb Z& \mbox{ for } k
= 0 ~\emph{\emph{or}}~n-1 ,\cr
\end{cases}  $
$H_{k}(A \cap L) = \begin{cases}0&
\mbox{for } \emph{\emph{others}}; \cr \mathbb Z& \mbox{ for } k = 0 ~\emph{\emph{or}}~n-2 .\cr
\end{cases}  $
\medskip

$ L $ is homeomorphic to $ \CC^{n-1} $, so $ H_{k}(L) = 0 $ for $ k\geq 1 $.

\bigskip

 Third, we have the Mayer-Vietoris sequence for $ A$ and  $L$:
\medskip

 $... \rightarrow H_{n-1}(A \cap L)
\rightarrow H_{n-1}(A) \oplus  H_{n-1}( L)  \rightarrow
H_{n-1}(A\cup L)
\rightarrow H_{n-2}(A \cap L)
\rightarrow H_{n-2}(A) \oplus  H_{n-2}( L)  \rightarrow...$
\medskip
which in the case $n>2$ is as follows:
$... \rightarrow 0
\rightarrow \Z \oplus 0  \rightarrow
H_{n-1}(A\cup L)
\rightarrow \Z
\rightarrow 0 \oplus  0  \rightarrow...$
\medskip

Therefore $H_{n-1}(A\cup L)   \cong\Z^{2}$.

The case $n=2$ is obvious.

The same arguments with $n$ replaced by any other dimension show that all
groups $\tilde H_i(A \cup L)$ with $i \neq n-1$ are trivial. \qed

\section{The fundamental group of the space of generic hyperplanes}

In this section we calculate the fundamental group $
\pi_1(\overset{\vee}{\CP^{n}}\setminus (\overset{\vee}{A} \cup
\overset{\ast}{A})) $, and in the next one we describe its action on $
H_{n-1}(A\cup L) $.
\medskip

\begin{theorem} 
\label{mainth}
If $n \geq 3$ then the group $\pi_1(\overset{\vee}{\CP^{n}}\setminus
(\overset{\vee}{A} \cup \overset{\ast}{A})$ is generated by three elements
$\alpha, \beta, \kappa$ with relations $\kappa \alpha =\beta \kappa$,
$\kappa^2=1$.
\end{theorem}

\begin{remark} \rm
Obviously, this presentation of the group can be reduced to one with only two
generators $\alpha, \kappa$ with the single relation $\kappa^2=1$. However, the
previous more symmetric presentation is more convenient for us.
\end{remark}

Denote by $  P\overset{\ast}{A}$ the set of all hyperplanes in
$\CP^{n-1}_{\infty},$ which are tangent to the hypersurface $\partial
 \overline{A}  \equiv \overline{A}\setminus A $ of ``infinitely distant'' points
 of $\overline{A}$.
\smallskip

Thus $ P\overset{\ast}{A} =  \overset{\vee}{(\overline{A}\setminus A)} $.
\medskip

Associating with any affine hyperplane in $\CC^n$ its infinitely distant part,
we obtain the down-left arrow in the commutative diagram of maps:

\begin{equation}\label{comdi} 
\begin{tikzcd}[column sep=6pc]
\overset{\vee}{\CP^{n}}\setminus (\overset{\vee}{A} \cup \overset{\ast}{A})
\arrow{r}{\emph{inclusion}}\arrow{d}{\CC^{1} \setminus\lbrace 2\emph{\emph{
points}}\rbrace} &
 \overset{\vee}{\CP^{n}}\setminus \overset{\ast}{A} \arrow{dl}{\CC^{1}} \\
 \overset{\vee}{\CP^{n-1}_\infty}\setminus P\overset{\ast}{A}
\end{tikzcd}
\end{equation}

Indeed, an affine hyperplane belongs to $\overset{\ast}{A}$ if and only if its
image under this map belongs to $ P\overset{\ast}{A}$. On the other hand, the
fiber of this map over any point of $\overset{\vee}{\CP^{n-1}_\infty}$ consists
of a pencil of affine hyperplanes parallel to one another, so it is a line
bundle. Any such fiber $\CC^1$ intersects the set $\overset{\vee}{A}$ at
exactly two points: indeed, for any non-asymptotic hyperplane there are exactly
two hyperplanes parallel to it and tangent to $A$.

Considering the fiber bundle represented by the left-hand part of the diagram
(\ref{comdi}),

\begin{eqnarray*}
\label{fifi}
\begin{CD}
E \ \ \ \\
\downarrow F \\
B \ \ \,
\end{CD}
\end{eqnarray*}
\smallskip
let $F= \CC^{1} \setminus\lbrace 2\emph{\emph{ points}}\rbrace, E=
\overset{\vee}{\CP^{n}}\setminus (\overset{\vee}{A} \cup \overset{\ast}{A}), B=
\overset{\vee}{\CP^{n-1}}\setminus P\overset{\ast}{A}$.

We have the exact sequence for the fiber bundle.
\begin{equation}
\label{exex} ... \rightarrow \pi_2(E) \rightarrow \pi_2(B) \rightarrow\pi_1(F)
\rightarrow \pi_1(E) \rightarrow \pi_1(B) \rightarrow \pi_0(F)
...\end{equation} \medskip

$ F $ is connected, so the rightmost arrow is trivial.

\begin{lemma} If $n >2$ then $\pi_1(B) = \Z_{2}$; if $n=2$ then $\pi_1(B)=\Z$.
\end{lemma}

Proof. The statement for $n=2$  is obvious: in this case $B$ is the complex
projective line less two points. For $n=3$ this statement follows by the
Zariski theorem (using the case $n=2$ as the base), see e.g. \cite{prasolov},
Chapter 6, \S 3. Finally, for $n>3$ it follows from the case $n=3$ by the strong
Lefschetz theorem, see \cite{GM}. \qed
\bigskip

\begin{lemma} Let $C$ be a smooth quadratic hypersurface in $\CP^{n-1}$.
If $n \neq 3$ then $\pi_2(\CP^{n-1} \setminus C)$ is trivial. $\pi_2(\CP^2 \setminus C) \cong \Z$. 
\end{lemma}

In particular, this is true for the base of our fiber bundle (\ref{fifi}).
\medskip

{\it Proof}. Let $[C] \subset \CC^n$ be the union of lines corresponding to the
points of $C$. We have a fiber bundle

\begin{eqnarray*}
\CC^{n}\setminus [C] \\
\downarrow \CC^{\ast} \\
\CP^{n-1}\setminus C
\end{eqnarray*}
\bigskip

 This fiber bundle is trivial because it is a restriction of the tautological
bundle of $\CP^{n-1}$ on the complement of a non-trivial divisor, so its first
Chern class is equal to 0.

Therefore $ \pi_2( \CC^{n}\setminus [C]) = \pi_2(\CP^{n-1}\setminus C) \oplus
\pi_2(\CC^{\ast}) = \pi_2(\CP^{n-1}\setminus C).$
\bigskip

Let $\varphi: \CC^n \to \CC$ be the quadratic polynomial defining the sets
$[C]$ and $C$. It defines the Milnor fibration $ \varphi: \CC^{n}\setminus
[C]\rightarrow \CC^{\ast}$.

Let $ E^{\prime}= \CC^{n}\setminus [C], B^{\prime}= \CC^{\ast}, F^{\prime}=
V_{\lambda}$. In this notation,  $ \pi_2(B)= \pi_2(\CP^{n-1}\setminus
C)=\pi_2(\CC^{n}\setminus [C]) =  \pi_2(E^{\prime})  $.

\medskip

We have the exact sequence for the fiber bundle.

$$...\pi_3(B^{\prime})\rightarrow\pi_2(F^{\prime})
\rightarrow
\pi_2(E^{\prime})
\rightarrow
\pi_2(B^{\prime}) \rightarrow\pi_1(F^{\prime})
\rightarrow
\pi_1(E^{\prime})
\rightarrow
\pi_1(B^{\prime})...$$
\bigskip

The base $B^\prime$ is homotopy equivalent to $S^1$, in particular the groups
$\pi_3(B^\prime)$ and $\pi_2(B^\prime)$ are trivial.
\medskip

Also, according to the Milnor theorem, $ F^{\prime} $ is homotopy equivalent to
$ S^{n-1} $.
\bigskip

Thus $ \pi_2(B)= \pi_2(E^{\prime})=\pi_2(F^{\prime}) = \pi_2( S^{n-1} ) = \begin{cases}0&
\mbox{ for } n \ne 3; \cr \mathbb Z& \mbox{ for } n=3.\cr
\end{cases} $
\qed

So for $n \neq 3$ the interesting fragment of the exact sequence
(\ref{exex}) reduces to

\begin{equation} \label{eee}
1 \rightarrow\pi_1(F)
\rightarrow
\pi_1(E)
\rightarrow
\pi_1(B) \rightarrow
1 .
\end{equation}
\bigskip

\begin{lemma} \label{excep}
In the case $n=3$ the map $\pi_2(E) \to \pi_2(B)$ in $($\ref{exex}$)$ is epimorphic.
\end{lemma}

{\it Proof.} By the construction of the generator of the group $\pi_2(B) \sim \Z$ in this case, this generator can be realised by the sphere consisting of complexifications of all oriented planes through the origin in $\R^3$. All these planes do not meet  the set $\overset{\vee}{A} \cup
\overset{\ast}{A}$, and hence define a 2-spheroid in $E$. \qed 
\medskip

So, the map $\pi_2(B) \to \pi_1(F)$  in (\ref{exex}) is trivial, and we can use the exact sequence (\ref{eee}) also in the case $n=3$.
\bigskip

$ \pi_1(F) = \Z\ast \Z,
\pi_1(B) = \Z_{2}$
\medskip

Thus $ \pi_1(E) $ has three generators $ \alpha, \beta, \kappa $, where
$\alpha$ and $\beta$ are two free generators of $\pi_1(F)$, and $\kappa$ is an
element of the coset $\pi_1(E) \setminus \pi_1(F)$.

We can realize these elements as follows. Choose the linear coordinates in
$\CC^n$ in which $A$ is given by the equation $z_1^2 + \dots + z_n^2 =1$. Take
for the base point in $\overset{\vee}{\CP^{n}}\setminus (\overset{\vee}{A} \cup
\overset{\ast}{A})$ the hyperplane $\{z_1=0\}$. The fiber $F$ containing this
point consists of all complex hyperplanes $\{z_1 = \mbox{const}\}$ parallel to
this one, they are characterized by the corresponding value of $z_1$. The
exceptional points of intersection with $\overset{\vee}{A}$ in this fiber
correspond to the values $1$ and $-1$.

Then for $\alpha$ and $\beta$ we take the classes of two simplest loops in
$\CC^1$ going along line intervals from 0 to  the points $1-\varepsilon$ (respectively,
$-1+\varepsilon$), $\varepsilon>0$ very small,  then turning counterclockwise around the point $1$ (respectively, $-1$) along a circle of radius $\varepsilon$, and coming back to
$0$.
\bigskip

For $\kappa$ we take the 1-parameter family of planes given by the equation
$(\cos \tau) z_1 + (\sin \tau) z_2=0$, $\tau \in [0,\pi]$.
\bigskip

\begin{lemma}
The element $\kappa$ thus defined does not belong to the image of $\pi_1(F)$ in
$\pi_1(E)$ under the second map in (\ref{eee}), i.e. its further map to $\pi_1(B)$ defines a generator of the latter group.
\end{lemma}

Indeed, it is easy to check this in the case $n=2$, which provides (via the
Zariski theorem) the generator of the latter group. \qed
\bigskip

The loop $\kappa$ defines also a loop in the base of our fiber bundle. Moving
the fibers over it and watching the corresponding movement of two exceptional
points, we get that $\kappa$ acts on $\pi_1(F)$ by permuting $\alpha$ and
$\beta$.

Theorem \ref{mainth} is proved. 

\section{Monodromy representation}

We know that
\begin{equation} \label{rh}
H_{n}(\CC^{n}, A\cup L) = \Z^{2},
\end{equation}
 see Proposition \ref{relativehomology}.

\begin{proposition}
For any $n$, the monodromy action of the group
$\pi_1(\overset{\vee}{\CP^{n}}\setminus (\overset{\vee}{A} \cup
\overset{\ast}{A}))$  on $H_n(\CC^n, A \cup L)$ has a 1-dimensional invariant
subspace.
\end{proposition}

{\it Proof.} This subspace is the image of the group $H_n(\CC^n,A) \cong \Z$ under the
obvious map $H_n(\CC^n, A) \to H_n(\CC^n, A \cup L)$; it corresponds via the boundary
isomorphism $H_n(\CC^n, A \cup L)\to H_{n-1}(A \cup L)$ in  (\ref{bpr}) to the image of the map $H_{n-1}(A) \to H_{n-1}(A \cup L)$. Indeed, this image does not depend on $L$.
\qed
\bigskip

It is convenient to fix the generators of this group (\ref{rh}) as follows.
Suppose again that $A$ is given by the equation
\begin{equation} \label{ball}
z_1^2 + \dots + z_n^2 =1,\end{equation} and the basepoint $L_0$ in the space of
planes is given by $z_1=0$. Then we have two relative cycles in $\CC^n$ (and
even in $\R^n$) modulo $A \cup L_0$: they are given by the two half-balls
bounded by the the surface ({\ref{ball}) and (the real part of) the hyperplane
$L_0$; we supply these half-balls with the orientations induced from a fixed
orientation of $\R^n$. It follows immediately from the proof of Proposition
\ref{relativehomology} that these two chains indeed generate the group
$H_{n}(\CC^{n}, A\cup L) = \Z^{2}.  $

Denote these two generators by $ a $ and $ b $. Namely, $a$ (respectively, $b$)
is the part placed in the half-space where $z_1>0$ (respectively, $z_1<0$). The
invariant subspace of the monodromy action is then generated by the sum of
these two elements: indeed, it is a relative cycle mod $A$ only.
\medskip

Let we study the action of loops $ \alpha, \beta,$ and $ \kappa $ on $ a $ and
$ b $.
\bigskip

\begin{proposition}
For any $n$, $ \kappa(a) = b, \kappa(b) = a.$
\end{proposition}

Proof. This follows immediately from the construction of both cycles $a$ and
$b$ and of the loop $\kappa$: when the hyperplane $L_\tau$ moves along this
loop, the parts of the space $\R^n$ bounded by the sphere (\ref{ball}) and real
parts of these hyperplanes move correspondingly and permute at the end of this
movement. \qed
\bigskip

\begin{proposition}
If $ n $ is odd, then the action of both loops $\alpha$ and $\beta$ is trivial.

If $ n $ is even, then
 $ \alpha(a) = -a,   \alpha(b) = 2a+b$, $ \beta(b) = -b,   \beta(a) = 2b+a$.
\end{proposition}

Proof. Both these statements follow immediately from the Picard--Lefschetz
formula, see Chapter III in \cite{APLT}. \qed
\bigskip

So, in the case of odd $n$ the monodromy action reduces to that of the group
$\Z_2$. In the case of even $n$ the monodromy group is infinite: for instance
the orbit of any generating element $a$ or $b$ consists of all points of the
integer lattice $\Z^2$ satisfying the conditions $u-v =1$ or $u-v=-1$.

\vspace{9pt}
Daodao Yang, Faculty of Mathematics, Higher School of Economics, Russia
\vspace{9pt}
\\
Address:  33/1 Studencheskaya, 121165, Moscow
\vspace{9pt}
\\
\textit{E-mail}: \textbf{dyang@edu.hse.ru}
\vspace{20pt}
 \\

\end{document}